\documentclass[10pt,reqno]{amsart}
\usepackage{amsmath}
\usepackage{amssymb}
\usepackage{amsthm}
\usepackage{eepic,epic}

\newlength{\defbaselineskip}
\newcommand{\setlinespacing}[1]%
           {\setlength{\baselineskip}{#1 \defbaselineskip}}

\numberwithin{equation}{section}

\newtheorem{thm}{Theorem}[section]

\newtheorem{cor}[thm]{Corollary}
\theoremstyle{definition}

\theoremstyle{remark}

\numberwithin{equation}{section}

\begin{document}

\title[Minimal support property]
{On minimal support properties of solutions of Schr\"odinger equations}

\author{Ihyeok Seo}


\address{School of Mathematics, Korea Institute for Advanced Study, Seoul 130-722, Republic of Korea}
\email{ihseo@kias.re.kr}

\maketitle

\begin{abstract}
In this paper we obtain minimal support properties of solutions of Schr\"odinger equations.
We improve previously known conditions on the potential for which the measure of
the support of solutions cannot be too small.
We also use these properties to obtain some new results on unique continuation
for the Schr\"odinger operator.
\end{abstract}

\section{Introduction}\label{sec1}

The main aim of this paper is to obtain the minimal support property
which implies a weak unique continuation for solutions of the stationary Schr\"odinger equation
\begin{equation*}
\Delta u(x)=V(x)u(x),
\end{equation*}
where $V$ is a potential function on $\mathbb{R}^n$, $n\geq3$.
The key ingredient in our method is
the following weighted $L^2$ inequality which is sometimes referred to as the Fefferman-Phong inequality:
For $u\in W^{1,2}(\mathbb{R}^n)$
\begin{equation}\label{ineq}
\int_{\mathbb{R}^n}|u|^2|V|dx\leq C(V)\int_{\mathbb{R}^n}|\nabla u|^2dx
\end{equation}
with a constant $C(V)$ depending on the potential $V(x)$.
Here, $W^{m,p}(\mathbb{R}^n)$ denotes the Sobolev space of functions
whose derivatives up to order $m$ belong to $L^p(\mathbb{R}^n)$.

Before stating our results, we shall recall some known results for \eqref{ineq} which
has played an important role in the study of the Schr\"odinger operator $-\Delta+V(x)$.
As it is well-known from \cite{F}, the inequality \eqref{ineq} holds for $V$ in the Fefferman-Phong class $\mathcal{F}^p$ for $1<p\leq n/2$.
This class was introduced by C. Fefferman and D. H. Phong to study eigenvalue estimates for the Schr\"odinger operator,
and is defined for $V\in L_{\textrm{loc}}^p$, $1\leq p\leq n/2$, by
$$V\in\mathcal{F}^p\quad\Leftrightarrow\quad
\|V\|_{\mathcal{F}^p}=\sup_{x,r}r^{2-n/p}\bigg(\int_{B(x,r)}|V(y)|^pdy\bigg)^{1/p}<\infty,$$
where $B(x,r)$ is the ball of radius $r>0$ centered at $x\in\mathbb{R}^n$.
In particular, $L^{n/2}=\mathcal{F}^{n/2}$ and $L^{n/2,\infty}\subset \mathcal{F}^p$ for all $1\leq p<n/2$.
So, if $p=n/2$, \eqref{ineq} follows easily from H\"older's inequality
and the Sobolev embedding theorem, with $C(V)=C\|V\|_{L^{n/2}}$.
On the other hand, it was noted in \cite{D} that \eqref{ineq} is not valid for $p=1$.
It was also shown in \cite{F} that for all $p>1$
$$\|V\|_{\mathcal{F}^1}\leq\|V\|\leq C\|V\|_{\mathcal{F}^p},$$
where $\|V\|$ denotes the least constant $C(V)$ for which ~\eqref{ineq} holds.
Later in \cite{KS}, Kerman and Sawyer showed that $\|V\|$ may be taken to be a constant multiple of
\begin{equation}\label{ksc}
\sup_{Q}\bigg(\int_Q|V(x)|dx\bigg)^{-1}\int_Q\int_Q\frac{|V(x)V(y)|}{|x-y|^{n-2}}dxdy<\infty.
\end{equation}
Here the sup is taken over all dyadic cubes $Q$ in $\mathbb{R}^n$, $n\geq3$.

Let us now define a wider class of potentials $V$ which allows \eqref{ineq}.

\medskip

\noindent\textbf{Definition.}
\textit{We say that $V$ is in the Kerman-Sawyer class $\mathcal{KS}_n$ if \eqref{ksc} holds.
Also we denote by $\|V\|_{\mathcal{KS}_n}$ the quantity in \eqref{ksc} and call it
the Kerman-Sawyer norm.
}

\subsection{Minimal support property}

Here we are mainly concerned with the minimal size of the support of solutions $u$
to the stationary Schr\"odinger equation
\begin{equation}\label{sse2}
\Delta u(x)=V(x)u(x),
\end{equation}
where $V$ is a potential function on $\mathbb{R}^n$, $n\geq3$.
In \cite{CEHL}, it was shown that if $u\in W_0^{1,2}(D)$ is a nontrivial solution of \eqref{sse2}
in a bounded domain $D\subset\mathbb{R}^n$,
then there exists a constant $C>0$ independent of $D$ such that
\begin{equation}\label{22}
|D|^{\frac2n-\frac1p}\|V\|_{L^p(D)}\geq C>0
\end{equation}
if $V\in L^p(D)$, $p>n/2$.
This implies that the measure $|D|$ cannot be too small.
(For an earlier result for the case $p=\infty$, see also \cite{CH}.)
The inequality \eqref{22} can be viewed as a relationship between a norm of the potential and
the measure of a domain containing the support of the solution.
This type of results is referred to as the \textit{minimal support property}.
Their method to obtain this property was based on
the Sobolev embedding $W_0^{1,2}(D)\hookrightarrow L^{2p'}(D)$,
where $p'$ is the H\"older conjugate exponent to $p$.
It was also showed in \cite{CEHL} that \eqref{22} does not hold
when $p<n/2$, in terms of counterexamples.

In this paper we extend these results to the class $\mathcal{KS}_n$ of potentials
which have locally small norms in the sense that
for all $z\in\mathbb{R}^n$
\begin{equation}\label{smallness}
\limsup_{r\rightarrow0}\|\textbf{1}_{B(z,r)}V\|_{\mathcal{KS}_n}\leq\varepsilon,
\end{equation}
where $\textbf{1}_A$ denotes the characteristic function of a set $A\subset\mathbb{R}^n$,
and $\varepsilon>0$ is a sufficiently small constant.
Note that the condition \eqref{smallness} is trivially satisfied for
$L^p$ norms, and recall that $L^{n/2}=\mathcal{F}^{n/2}\subset\mathcal{KS}_n$.
Hence, this makes it possible to improve \eqref{22} to the endpoint case $p=n/2$.
Also, we consider the problem for a more general differential inequality
\begin{equation}\label{ssi0}
|\Delta u(x)|\leq|V(x)u(x)|.
\end{equation}
The method here will be based on ~\eqref{ineq} which can be viewed as the weighted embedding
$W^{1,2}(\mathbb{R}^n)\hookrightarrow L^{2}(V)$,
and can be also applied to the magnetic Schr\"odinger operator (see Section \ref{sec3}).

Our first result is the following.

\begin{thm}\label{thm}
Let $u\in W^{2,2}(\mathbb{R}^n)$ be a solution of \eqref{ssi0},
whose support is contained in a ball $B(a,r)$.
If $\textbf{1}_{B(a,r)}V\in\mathcal{KS}_n$,
then there exists $r_0>0$ such that for all $r<r_0$
$$\|\textbf{1}_{B(a,r)}V\|_{\mathcal{KS}_n}\geq C>0$$
with a constant $C$ independent of $r$.
\end{thm}

As an immediate consequence of this theorem, we have the following
minimal support property.

\begin{cor}\label{cor}
Let $u\in W^{2,2}(\mathbb{R}^n)$ be a solution of \eqref{ssi0}
with $V\in\mathcal{KS}_n$ satisfying \eqref{smallness}.
If the support of $u$ is contained in a ball $B(a,r)$,
then the measure $|B(a,r)|$ cannot be too small.
\end{cor}

Let us now consider a different class of potentials, the Kato class $\mathcal{K}_n$,
which is defined for $V\in L_{\textrm{loc}}^1$ by
\begin{equation*}
V\in \mathcal{K}_n\quad\Leftrightarrow\quad\lim_{r\rightarrow0}
\sup_{x\in\mathbb{R}^n}\int_{|x-y|<r}\frac{|V(y)|}{|x-y|^{n-2}}dy=0.
\end{equation*}
(This class named for T. Kato \cite{K2} has arisen in the analysis of self-adjointness
of the Schr\"odinger operator.)
We denote
\begin{equation*}
\eta(r,V)\equiv\sup_{x\in\mathbb{R}^n}\int_{|x-y|<r}\frac{|V(y)|}{|x-y|^{n-2}}dy.
\end{equation*}
Then, a variant of the Kato class denoted by $\widetilde{\mathcal{K}}_n$ can be defined as
\begin{equation*}
V\in \widetilde{\mathcal{K}}_n\quad\Leftrightarrow\quad\eta(r,V)<\infty\quad\text{for all}\,\, r>0,
\end{equation*}
and it is possible to obtain ~\eqref{ineq} with $C(V)=C\eta(2r,V)$ for this class (see \cite{Z}):
For $u\in C_0^\infty(\mathbb{R}^n)$ supported in $B(z,r)$,
\begin{equation}\label{ineq2}
\int_{\mathbb{R}^n}|u|^2|V|dx\leq C\eta(2r,V)\int_{\mathbb{R}^n}|\nabla u|^2dx.
\end{equation}
Making use of \eqref{ineq2}, we obtain the following result,
which can be also seen as extensions to the class $\widetilde{\mathcal{K}}_n$
of the previous results (see \eqref{22}) in \cite{CH,CEHL}
because $L^p\subset\mathcal{K}_n\subset\widetilde{\mathcal{K}}_n$ for all $p>n/2$.

\begin{thm}\label{thm2}
Let $u\in W^{2,2}(\mathbb{R}^n)$ be a solution of \eqref{ssi0},
whose support is contained in a ball $B(a,r)$.
If $\textbf{1}_{B(a,r)}V\in\widetilde{\mathcal{K}}_n$,
then there exists $r_0>0$ such that for all $r<r_0$
$$\eta(2r,V)\geq C>0$$
with a constant $C$ independent of $r$.
\end{thm}

Note that $\lim_{r\rightarrow0}\eta(r,V)=0$ for $V\in \mathcal{K}_n$.
Hence, as in Corollary \ref{cor}, Theorem \ref{thm2} directly implies the following.

\begin{cor}\label{cor2}
Let $u\in W^{2,2}(\mathbb{R}^n)$ be a solution of \eqref{ssi0} with $V\in\mathcal{K}_n$.
If the support of $u$ is contained in a ball $B(a,r)$,
then the measure $|B(a,r)|$ cannot be too small.
\end{cor}

It is clear that the same minimal support results as in Corollaries \ref{cor} and \ref{cor2}
hold for the eigenvalue equation of the Schr\"odinger operator,
\begin{equation}\label{eigen}
(-\Delta+V(x))u=Eu,
\end{equation}
because $V-E$ satisfies \eqref{smallness} if $V$ satisfies it (see \eqref{33})
and similarly $V-E\in K_n$ if $V\in K_n$.

\subsection{Applications to unique continuation}

Now we would like to emphasize that our minimal support results imply
some new results on unique continuation for the Schr\"odinger operator.

Given a partial differential operator $P(x,D)$, we say that it has the \textit{unique continuation property}
if the solution $u$ of $Pu=0$ which vanishes in an open subset of its domain of definition
must vanish identically.
Historically, the study of the unique continuation originated from its connection with the
uniqueness of the Cauchy problem, to which, in many cases, it is equivalent.
On the other hand, the most important motivation came from quantum mathematical physics.
For example, Kato \cite{K} proved that if $V$ has a compact support, then for $E>0$,
all the corresponding eigensolutions $u$ of \eqref{eigen} must vanish
outside of a ball. Hence, the unique continuation implies a proof of absence
of positive eigenvalues $E$.
At this point, it is enough to establish a \textit{weak unique continuation property}
which says that the solution that vanishes in the complement of a compact subset must vanish identically.
Of course, the unique continuation directly implies the weak one.

The first result on the unique continuation for the Schr\"odinger operator
is due to Carleman \cite{C}, who showed it in dimension two
for bounded potentials $V\in L_{\textrm{loc}}^\infty(\mathbb{R}^2)$.
This was extended to higher dimensions $n\geq3$ by M\"uller \cite{M}.
Since then, a great deal of work was devoted to the unbounded cases
$V\in L_{\textrm{loc}}^p(\mathbb{R}^n)$, $p<\infty$.
Among others, Jerison and Kenig \cite{JK} proved
the unique continuation for the differential inequality
\begin{equation}\label{wuc}
|\Delta u(x)|\leq|V(x)u(x)|
\end{equation}
if $V\in L_{\textrm{loc}}^{n/2}(\mathbb{R}^n)$ for $n\geq3$,
and $V\in L_{\textrm{loc}}^p(\mathbb{R}^2)$ for $p>1$.
This result later turns out to be optimal in the context of $L^p$ potentials (\cite{KN,KT}).
In this regard, the later developments have been made
to extend these $L^p$ potentials to more singular ones such as Kato and Fefferman-Phong potentials.
In fact, the smallness condition like \eqref{smallness} has already appeared in these developments.
In~\cite{St}, Stein showed the unique continuation for ~\eqref{wuc} if
$V\in L_{\textrm{loc}}^{n/2,\infty}(\mathbb{R}^n)$, $n\geq3$, with the norm in \eqref{smallness}
replaced by $L^{n/2,\infty}$-norm.
This was extended by Chanillo and Sawyer \cite{CS}
to the Fefferman-Phong potentials $V\in\mathcal{F}^p$
for $p>(n-1)/2$ if $n\geq3$, and $p>1$ if $n=2$,
with the norm in \eqref{smallness} replaced by $\mathcal{F}^p$-norm.
(This was also improved by Wolff \cite{W} to $p>(n-2)/2$, $n\geq4$.)
In low dimensions $n\leq3$, the unique continuation was further extended by them
to the class $\mathcal{KS}_n$ with \eqref{smallness},
and was obtained by Sawyer \cite{S} for the Kato class $\mathcal{K}_n$.

On the other hand, there have been no results on the weak unique continuation for~\eqref{wuc}
beyond trivial ones that follows directly from the above unique continuation results.
As an immediate consequence of Corollaries \ref{cor} and \ref{cor2},
we have the following new result on the weak unique continuation.

\begin{cor}
Let $u\in W^{2,2}(\mathbb{R}^n)$ be a solution of \eqref{wuc}
which vanishes in the complement of a compact set $K$.
If $K$ is contained in a sufficiently small ball,
then $u$ must vanish identically provided that
$V\in\mathcal{KS}_n$ with \eqref{smallness}, or $V\in\mathcal{K}_n$.
\end{cor}

Throughout this paper, the letter $C$ stands for constants possibly different at each occurrence.


\section{Proof of Theorems \ref{thm} and \ref{thm2}}

The method of proof of Theorems \ref{thm} and \ref{thm2} follows the same argument,
which is based on the weighted embeddings \eqref{ineq} and \eqref{ineq2}, respectively.
So, we omit the proof for Theorem \ref{thm2}.

Now we prove Theorem \ref{thm}.
Since $u\in W^{2,2}(\mathbb{R}^n)$,
we can find functions $u_m\in C_0^\infty(\mathbb{R}^n)$ so that if $m\rightarrow\infty$
\begin{equation}\label{conv}
u_m\rightarrow u\quad\text{in}\quad W^{2,2}(\mathbb{R}^n).
\end{equation}
Now, let $\psi$ be a smooth cut-off function such that $\psi=1$ on $B(a,r)$
and $\psi=0$ on $\mathbb{R}^n\setminus B(a,2r)$.
Then it follows from \eqref{ineq} that
\begin{align}\label{55}
\nonumber\int_{B(a,r)}|u_m|^2|V|dx&=\int|u_m\psi|^2|\textbf{1}_{B(a,r)}V|dx\\
&\leq C\|\textbf{1}_{B(a,r)}V\|_{\mathcal{KS}_n}\int|\nabla(u_m\psi)|^2dx.
\end{align}
By using Green's identity, we also see that
\begin{align}\label{333}
\nonumber\int|\nabla(u_m\psi)|^2dx
\nonumber&=\int_{B(a,2r)}\nabla(u_m\psi)\nabla(\overline{u_m\psi})dx\\
&=-\int_{B(a,2r)}\Delta(u_m\psi)\overline{u_m\psi}dx.
\end{align}
Then, we write
\begin{align}\label{000}
\int_{B(a,2r)}\Delta(u_m\psi)\overline{u_m\psi}dx=&
\int_{B(a,2r)}\overline{u_m\psi}\psi\Delta u_mdx\\
\nonumber&+\int_{B(a,2r)}2\overline{u_m\psi}\nabla u_m\cdot\nabla\psi dx
+\int_{B(a,2r)}\overline{u_m\psi}u_m\Delta\psi dx.
\end{align}
For the first term on the right-hand side of \eqref{000}, we note that
\begin{equation}\label{357}
\lim_{m\rightarrow\infty}\int_{B(a,2r)}\overline{u_m\psi}\psi\Delta u_mdx
=\int_{B(a,2r)}\overline{u\psi}\psi\Delta u dx.
\end{equation}
Indeed, from \eqref{conv} and the simple fact that
\begin{align*}
\int_{B(a,2r)}\overline{u_m\psi}\psi\Delta u_m
-\overline{u\psi}\psi\Delta u dx
\leq &C\|u_m-u\|_{L^2}\|\Delta u_m\|_{L^2}\\
&+C\|u\|_{L^2}\|\Delta u_m-\Delta u\|_{L^2},
\end{align*}
we get \eqref{357}.
Similarly, by letting $m\rightarrow\infty$ and using the support properties of $u$, $\psi$,
the last two terms on the right-hand side of \eqref{000} become zero.
Consequently,
\begin{align*}
\lim_{m\rightarrow\infty}\int_{B(a,2r)}\Delta(u_m\psi)\overline{u_m\psi}dx
&=\int_{B(a,2r)}\overline{u\psi}\psi\Delta udx\\
&=\int_{B(a,r)}\overline{u}\Delta udx,
\end{align*}
since we are assuming $\text{supp}\,u\subset B(a,r)$, and $\psi=1$ on $B(a,r)$.
Combining \eqref{55}, \eqref{333} and this, we get
\begin{equation}\label{898}
\lim_{m\rightarrow\infty}\int_{B(a,r)}|u_m|^2|V|dx\leq
C\|\textbf{1}_{B(a,r)}V\|_{\mathcal{KS}_n}\int_{B(a,r)}|u||\Delta u|dx.
\end{equation}
Next, we note that
$$\bigg(\int_{B(a,r)}|u|^2|V|dx\bigg)^{1/2}
\leq\bigg(\int_{B(a,r)}|u-u_m|^2|V|dx\bigg)^{1/2}
+\bigg(\int_{B(a,r)}|u_m|^2|V|dx\bigg)^{1/2}.$$
Since the first term on the right-hand side is bounded by
$$\|\textbf{1}_{B(a,r)}V\|_{\mathcal{KS}_n}^{1/2}\bigg(\int_{B(a,r)}|\nabla(u-u_m)|^2dx\bigg)^{1/2},$$
from \eqref{conv} and \eqref{898}, we see that
\begin{equation}\label{111}
\int_{B(a,r)}|u|^2|V|dx\leq
C\|\textbf{1}_{B(a,r)}V\|_{\mathcal{KS}_n}\int_{B(a,r)}|u||\Delta u|dx.
\end{equation}
Using this and \eqref{ssi0}, we conclude that
\begin{equation}\label{88}
\int_{B(a,r)}|u|^2|V|dx\leq C\|\textbf{1}_{B(a,r)}V\|_{\mathcal{KS}_n}\int_{B(a,r)}|u|^2|V|dx.
\end{equation}

Now, we want to show that
\begin{equation}\label{33}
\lim_{r\rightarrow0}\|\textbf{1}_{B(a,r)}\|_{\mathcal{KS}_n}=0
\end{equation}
which implies that there exists $\widetilde{r}>0$ so that
$\textbf{1}_{B(a,r)}\in\mathcal{KS}_n$ for all $r<\widetilde{r}$.
To see this, we first note that
\begin{equation*}
\|V\|_{\mathcal{KS}_n}\leq
\sup_{x\in\mathbb{R}^n}\int_{\mathbb{R}^n}\frac{|V(y)|}{|x-y|^{n-2}}dy.
\end{equation*}
Hence, we only need to show
$$\lim_{r\rightarrow0}\sup_{x\in\mathbb{R}^n}\int_{|y-a|<r}\frac{1}{|x-y|^{n-2}}dy=0.$$
But, this is an easy computation by noting that
\begin{align*}
\sup_{x\in\mathbb{R}^n}\int_{|y-a|<r}\frac{1}{|x-y|^{n-2}}dy
&\leq\sup_{|x-a|<2r}\int_{|y-a|<r}|x-y|^{-(n-2)}dx\\
&\qquad\qquad\qquad+\sup_{|x-a|\geq2r}\int_{|y-a|<r}r^{-(n-2)}dx\\
&\leq\sup_{x\in\mathbb{R}^n}\int_{|x-y|<4r}|x-y|^{-(n-2)}dx+Cr^2\\
&\leq Cr^2.
\end{align*}
Then, since \eqref{ssi0} is also satisfied for the potential $\widetilde{V}=|V|+\textbf{1}_{B(a,r)}$,
from \eqref{88}, we see that for all $r<\widetilde{r}$
\begin{equation}\label{888}
\int_{B(a,r)}|u|^2\widetilde{V}dx\leq
C\|\textbf{1}_{B(a,r)}|V|+\textbf{1}_{B(a,r)}\|_{\mathcal{KS}_n}\int_{B(a,r)}|u|^2\widetilde{V}dx.
\end{equation}
By deleting the term
$$0<\int_{B(a,r)}|u|^2\widetilde{V}dx
\,\,\Big(\leq C\|\textbf{1}_{B(a,r)}\widetilde{V}\|_{\mathcal{KS}_n}\int_{B(a,r)}|\nabla u|^2dx\Big)\,<\infty$$
from \eqref{888}, we conclude that
$$\|\textbf{1}_{B(a,r)}V\|_{\mathcal{KS}_n}+\|\textbf{1}_{B(a,r)}\|_{\mathcal{KS}_n}\geq C>0.$$
Now, from this and \eqref{33}, there exists $r_0>0$ such that for all $r<r_0$
$$\|\textbf{1}_{B(a,r)}V\|_{\mathcal{KS}_n}\geq C/2>0.$$
This completes the proof.

\section{Concluding remarks}\label{sec3}

The method here can be also applied to obtain similar results for the magnetic Schr\"odinger operator
$-\Delta_{\vec{A}(x)}+V(x)$, where $\vec{A}(x)=(A_1(x),...,A_n(x))$ is a magnetic potential,
such that $A_j(x)$, $j=1,...,n$, are real valued functions,
and $\Delta_{\vec{A}(x)}$ denotes the magnetic Laplacian defined by
\begin{align*}
\Delta_{\vec{A}(x)}&=\sum_j(\partial_j+iA_j)^2\\
&=\Delta+2i\vec{A}\cdot\nabla+idiv\vec{A}-|\vec{A}|^2.
\end{align*}
Setting $\vec{B}=-2i\vec{A}$ and $W=-idiv\vec{A}+|\vec{A}|^2+V$, we can rewrite the operator as
$-\Delta+\vec{B}\cdot\nabla+W$.
Hence we reduce the problem to the differential inequality
\begin{equation}\label{mSch}
|\Delta u|\leq|Wu|+|\vec{B}\cdot\nabla u|,
\end{equation}
where $W:\mathbb{R}^n\rightarrow\mathbb{C}$ and $\vec{B}=(B_1,...,B_n):\mathbb{R}^n\rightarrow\mathbb{C}^n$.

Then, with the potential $V$ replaced by $|W|+|\textrm{div}\vec{B}|$,
all the results in Section \ref{sec1} can be valid for \eqref{mSch}.
Indeed, from integration by parts, note that
\begin{align*}
\int_{B(a,r)}2u\nabla u\cdot\vec{B}dx&=\int_{B(a,r)}\nabla(u^2)\cdot\vec{B}dx\\
&=-\int_{B(a,r)}u^2\,\textrm{div}\vec{B}dx.
\end{align*}
Then, the integral on the right-hand side of \eqref{111} can be handled by
\begin{align*}
\int_{B(a,r)}|u||\Delta u|dx
&\leq\int_{B(a,r)}|u|^2|W|+|u\nabla u\cdot\vec{B}|dx\\
&\leq C\int_{B(a,r)}|u|^2(|W|+|\textrm{div}\vec{B}|)dx.
\end{align*}
Hence, we get \eqref{88} for $V=|W|+|\textrm{div}\vec{B}|$.
Now, the remaining part follows easily from the same argument.
We omit the details.


\end{document}